\documentclass[a4paper,12pt,reqno]{amsart}

\usepackage{amsmath}%
\usepackage{amsfonts}%
\usepackage{amssymb}%
\usepackage{graphicx}

\theoremstyle{plain}

\newtheorem{conclusion}{Conclusion}

\newtheorem{criterion}{Criterion}

\newtheorem{lemma}{Lemma}

\newtheorem{problem}{Problem}

\newtheorem{remark}{Remark}

\numberwithin{equation}{section}

\begin{document}

\title[General Methods of search]{General method of \textit{Complex Polynomial} for determining the radius of the circle circumscribed to a cyclic polygon an arbitrary number of sides, and some important consequences.}
\author{\textit{Denis Mart\'{i}nez T\'{a}panes}}
\address[A. One and A. Two]
{University of Medical Sciences of Villa Clara, Cuba 
\newline%
\indent Central University of Las Villas, Cuba }%
\email[A. One]{denismt@ucm.vcl.sld.cu}%
\urladdr{http://www.ucm.vcl.sld.cu/}
\author{\textit{J. Enrique Mart\'{i}nez Serra}}
\email[A.~Two]{josee@uclv.edu.cu}%
\urladdr{http://www.uclv.edu.cu/}
\author{\textit{L. Osiel Rodr\'{i}guez Ca$\tilde{n}$izarez}}
\email[A.~Three]{lrcanizares@uclv.cu}%
\urladdr{http://www.uclv.edu.cu/}
\thanks{This paper is in final form and no version of it will be submitted for publication elsewhere.}
\date{June 28, 2015}
\dedicatory{\normalsize Dedicated to our Professor Eberto Morgado.}

\begin{abstract}

This paper presents a general method for obtaining radius of the corresponding circumference to a cyclical polygon $n$ sides given the lengths of said sides, using the notion of complex number. As of radius $r$, obtained, can then be calculated polygon area in question applying known expressions which requires, after $n=4$, a powerful calculation tool. Are also given like elements regarding non convex polygons cyclic, although in this respect it deepens less and is not considered, course, the calculation of areas.
\end{abstract}
\maketitle

\section{Introduction}

Determining the areas of the cyclic polygons from its sides is usually a very complex task and solving this is not given. By definition, any cyclic polygon is inscribed in a certain unique circumference and the problem of finding the area of the polygon from the lengths of the sides it begins by determining the radius of the circle in which the polygon is inscribed. Once you found this radius, area will be found by the sum of the areas of triangles having two sides equal to the radius and the other equal to the corresponding side of the polygon; it is evident that the polygon of $n$ sides it is constituted for $n$ triangles of this type in turn.

Taking advantage of the current possibilities of calculation using computers, calculations related to this problem do not present serious difficulties provided that exist the appropriate theoretical method to be implemented in a package like Mathematica; but as will be seen below, particular expressions of area in terms of sides of polygon passing of four sides largely lose interest because they can not give clear information due to its extreme length and complexity, unless you want to see a curious thing.

In this paper we simply report the method conceived by the authors and put enough examples of its effectiveness. This paper presents a general method for obtaining radius of the corresponding circumference to a cyclical polygon (convex ) n-sided given the lengths of said sides, using the notion of complex number. As of radius $r$, obtained, can then be calculated polygon area in question applying, for example, the expression $\frac{1}{4} \sum_{1}^{n} l_i \sqrt{4r^2-(l_i )^2 }$ (If the center of circumcircle stand in to its area) which requires, after $n=4$, a powerful calculation tool. The contrary case (The center of circumcircle not stand in to the polygon area) also it addressed with slight variations in the expressions required for the calculation.

The article consists of three general sections namely the first section it carry out the method for calculating the radius of the circumscribed circumference, the second deals with the calculation of areas based on the calculation of the radius, and the third has sufficient examples to demonstrate the accuracy of the methods set first.

We clarify that each of the specific examples that were given numerical values of the sides were checked with a ruler and compass in the case of calculating the radio and some practical physical method for the case of the areas. For this we have a certain area pattern which in this case was that of a sheet of paper from high quality ($216 \times 279$) to which the weight is measured; on sheets of the same type figures whose area was sought were outlined; then the figures were cut and weighed carefully with the same analytical balance with which the pattern is weighed. The results were always very accurate.

Are also given like elements regarding non convex polygons cyclic, although in this respect it deepens less and is not considered the calculation of areas. But the formula used to implement a method analogous to the convex case is established, as well as some general ideas of how to interpret the results of such implementation. At least it is clear that in case the search does not convex radius required to know certain specific characteristics of cyclic polygon and can not hope for a generic formula as in the case of convex. This means there is some general formula parameter that is unique to the convex, but it is not for the non-convex.

\section{Necessary theoretical development for calculating the radius of the circumscribed circle}

\subsection{Convex cyclic polygons:}

Suppose we have a circle and a convex polygon of n sides inscribed on it as it is shown in figure 1. The relationship between $\alpha_k$ and the corresponding side $l_k$ is (under the law of cosines) the following,

\begin{figure}[h]
	\centering
		\includegraphics[width=0.30 \textwidth]{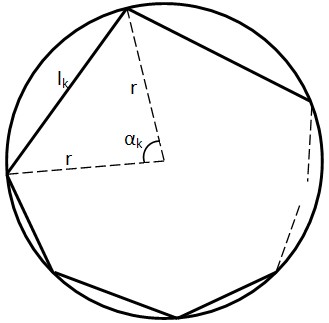}
	\caption{Partial picture of the polygon inscribed with the side and the corresponding angle.}
	\label{fig:1}
\end{figure}

\begin{equation}
\cos\alpha_k=1-\frac{l_k^{2} }{2r^2}
\label{eqn20.10}
\end{equation}

, from this it follows

\begin{equation}
\sin\alpha_k=\pm\sqrt{1-\left(1-\frac{l_k ^2}{2r^2}\right)^2}
\end{equation}

If the centre of the circle is located within the area bounded by the polygon (below a classification for this will be required) will be fulfilled that 

\begin{equation}
	\sum_1^n\alpha_k =2\pi
  \label{eqn20.30}
\end{equation}

It is equivalent to following,

\begin{equation}
	\alpha_n=2\pi-\sum_{1}^{n-1}\alpha_k
\end{equation}

While if the opposite occurs, and $l_n$ is not less than any of the other side's takes place as follows,

\begin{equation}
\alpha_n=\sum^{n-1}_{1}\alpha_k
\label{eqn20.50}
\end{equation}

However, both expressions are equivalent with respect to evaluation of trigonometric functions, which can be seen in the following identities,

\begin{equation}
\cos\alpha_n=\cos\left(2\pi-\sum_{1}^{n-1}\alpha_k\right)=\cos\left(\sum_{1}^{n-1}\alpha_k\right)=\cos\left(-\sum_{1}^{n-1}\alpha_k\right)
\label{eqn20,60}
\end{equation}

\begin{equation}
\sin\alpha_n=\sin\left(2\pi-\sum_{1}^{n-1}\alpha_k\right)=\sin\left(-\sum_{1}^{n-1}\alpha_k\right)
\label{eqn20,70}
\end{equation}

Consider now the following complex number (can to note that the module of the number given to the right of the following relationship is equal $1$, such and like it should be in this case).

\begin{equation}
	e^{i\alpha_k}=1-\frac{l_k ^2}{2r^2}+i\sqrt{1-\left(1-\frac{l_k ^2}{2r^2}\right)^2}
\end{equation}

Where can shape the next product:

\begin{equation}
\prod^{n-1}_{k=1}e^{i\alpha_k}=\prod^{n-1}_{1}\left(1-\frac{l_k ^2}{2r^2}+i\sqrt{1-\left(1-\frac{l_{k}^2}{2r^2}\right)^2 }\right)
\end{equation}

However it is clear that from (\ref{eqn20,60}) and (\ref{eqn20,70}) may be considered in general $\alpha_n=-\sum^{n-1}_{1}\alpha_k $ and therefore

\[\prod^{n-1}_{k=1}e^{i\alpha_k}=e^{i\sum^{n-1}_{1}\alpha_k}=e^{-i\alpha_n}\]

\begin{equation}
\prod^{n}_{k=1}e^{i\alpha_k}=e^{i\sum^{n}_{1}\alpha_k}=1
\end{equation}

Which ultimately allows us to propose the following,

\begin{equation}
	1-\prod^{n}_{k=1}\left(1-\frac{l_k ^2}{2r^2}+i\sqrt{1-\left(1-\frac{l_k ^2}{2r^2}\right)^2 }\right)=0 \label{eqn20.110}
\end{equation}

Taking the real part, or the complex, of above relation, it will be considering an equation, with radicals, which allows find the radius of the circle where you can join the polygon of $n$ sides considered.

\begin{remark}
Whenever using this method (When the radius is searched alone for the convex case), must be despised all solutions (Subsection 2.2 show that they correspond to all the possible not convex cases) except one of them which will be neither negative nor complex; and between the positive real roots should be chosen within the set $\Phi=[\left\{r_i \right\},\,r_i>(\sum_1^n l_k )/2\pi].$ It would be extremely interesting to determine if the set $\Phi$ it contains one or more elements different from the empty set, although it is not objective of this article and it won't be tried here. 
\end{remark}

The previous remark the following open problem arises

\begin{problem}
Determine whether the number of elements in the set  $\Phi$ is always the same regardless of the number of sides. In that case, which is that amount?.
\end{problem}

\subsection{Generalization to non convex cyclic polygons:}

In the same way as was done in the case convex can be formed the following complex numbers,

\begin{equation}
N_{p}(r)=f_{p}(r)+i\sqrt{1-f^{2}_{p}(r)}=e^{\sum^{a_{p}}_{i=1}\alpha_{ip}}
\label{eqn20.120}
\end{equation}

were,

\begin{center}
\[f_{1}(r)=\cos\left(\sum^{a_{1}}_{i=1}\alpha_{i1}\right)=1-\frac{l^{2}_{1}}{2r^{2}},\ \sin\left(\sum^{a_{1}}_{i=1}\alpha_{i1}\right)=\sqrt{1-f^{2}_{1}(r)}\]
\[f_{2}(r)=\cos\left(\sum^{a_{2}}_{i=1}\alpha_{i2}\right)=1-\frac{l^{2}_{2}}{2r^{2}},\ \sin\left(\sum^{a_{2}}_{i=1}\alpha_{i1}\right)=\sqrt{1-f^{2}_{2}(r)}\]
\[\vdots\] 
\[f_{n}(r)=\cos\left(\sum^{a_{n}}_{i=1}\alpha_{in}\right)=1-\frac{l^{2}_{n}}{2r^{2}},\ \sin\left(\sum^{a_{n}}_{i=1}\alpha_{i1}\right)=\sqrt{1-f^{2}_{n}(r)}\] 
\end{center}

However in the general case every side of the polygon may comprise various angles (see fig.2) which is represented by the above expressions in the respective sums.

From (\ref{eqn20.120}) you are obviously obtained,

\begin{equation}
\prod^{n}_{p=1}N_{p}\left(r\right)=e^{i\left(\sum^{a_{1}}_{i=1}\alpha_{k_{i1}}+\sum^{a_{2}}_{i=1}\alpha_{k_{i2}}+\ldots+\sum^{a_{n}}_{i=1}\alpha_{k_{in}}\right)}=e^{i\left(q_{1}\alpha_{1}+q_{2}\alpha_{2}+...+q_{n}\alpha_{n}\right)}
\label{eqn20.130}
\end{equation}

Where $q_{i}$ are non-negative integers three of which at least can not be null. Moreover it has to,

\begin{equation}
e^{i\left(q_{1}\alpha_{1}+q_{2}\alpha_{2}+...+q_{n}\alpha_{n}\right)}=e^{iE2\pi}
\label{eqn20.140}
\end{equation}

Where $E$ is a positive real number, which can be displayed knowing that $\alpha_{i}=\delta_{i}2\pi$ were $\delta_{i}$ it is a real fraction of the unit. 

And ultimately that will be

\begin{equation}
\prod^{n}_{p=1}N_{p}(r)=e^{iE2\pi}
\label{eqn20.150}
\end{equation}

Notice how in the convex case $q_{i}$ are all equal to unity and the $E$ is obviously $1$ obtaining (\ref{eqn20.110}) immediately.

\begin{figure}[h]
\centering
\includegraphics[width=0.30\textwidth]{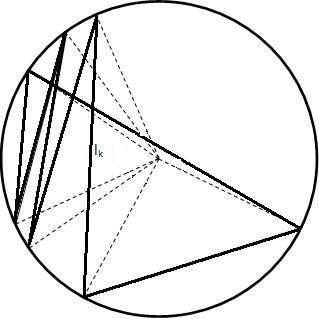}
\caption{Crossed seven-sided polygon. Note that the $l_{k}$
will correspond five successive angles.}
\label{fig:2}
\end{figure}

The following lemma will be important to address the non-convex case:

\begin{lemma}
The integer value of $q_{i}$ (See (\ref{eqn20.130})) is equal the number of sides of the polygon in question passing the sector of the circumference corresponding to the angle $\alpha_{i}$.
\end{lemma}

\begin{proof}
It is only necessary to refer to (\ref{eqn20.120}) and note that the number of times $\alpha_{i}$ it is repeated in the exponent of (\ref{eqn20.130}), which it is precisely $q_{i}$, it equals the number of equalities of the form (\ref{eqn20.120}) in that the angle appears; however it is clear that if $\alpha_{i}$ it appears in a quantity $m$ of these equalities it is that an amount of $ m $ sides pass through its sector by virtue of the way $f_{p}(r)$ it has . The lemma has been demonstrated.
\end{proof}

\subsection{Implementation in a mathematical package (Convex Case):} 
The following procedure exemplified in eight different cases (In section $4$), where $r$ represents the radius that is sought, it was implemented in the math package Mathematica $9.0$:

\begin{enumerate}
\item Since the polygon $n$ sides $l_1, l_2,...,l_n$; do the following

\[Expand[1-\prod^{n}_{k=1}\left(1-\frac{l_k^2}{2r^2}+i\sqrt{1-\left(1-\frac{l_k ^2}{2r^2}\right)^2}\right)]\]

The result, it will be called $C(r,l_{1},l_{2},...,l_{n})$, it will be a ``disorganized'' complex number.
\item Perform Re[$C(r,l_{1},l_{2},...,l_{n})$].
The result is a complex number, call it $A(r,l_{1},l_{2},...,l_{n})+i B(r,l_{1},l_{2},...,l_{n})$, where the real and the imaginary part are well defined. However we not free ourselves from $\emph{i}$ because the above command does not define the real character of the elements $r,l_{1},l_{2},...,l_{n}$.
\item Take $A(r,l_{1},l_{2},...,l_{n})$ and perform, according to the shape of the desired result, one of the commands

Solve[$A(r,l_{1},l_{2},...,l_{n}),r$]; NSolve[$A(r,l_{1},l_{2},...,l_{n}),r$]

\item Choosing (For the specific convex case) the right root, discarding negative, complex and that, multiplied by $2\pi$, they give a numerical result smaller or equal to the perimeter of the polygon. This is based on the fact that the perimeter of the circumscribed circle has to be larger than the inscribed polygon that does not exceed the border.
\end{enumerate}

\begin{remark}
Of course the roots generated by the above procedure will not only to the convex polygons, but the corresponding to crossed polygons for which $q_{i}$ are all equal (see (\ref{eqn20.130})); example of regular stars five- and seven points will be shown.
\end{remark}

\begin{remark}
For regular convex polygons in the expression given point $ 1 $ of the previous implementation can be simplified significantly by the following transformation,

\[\left\{\prod^{n}_{k=1}\left(1-\frac{l^2}{2r^2}+i\sqrt{1-\left(1-\frac{l^2}{2r^2}\right)^2}\right)=1\right\}\equiv\]
\[\equiv\left\{\left(1-\frac{l^2}{2r^2}+i\sqrt{1-\left(1-\frac{l^2}{2r^2}\right)^2}\right)^{n}=e^{i2\pi}\right\}\equiv\]
\[\equiv\left\{1-\frac{l^2}{2r^2}+i\sqrt{1-\left(1-\frac{l^2}{2r^2}\right)^2}=e^{i\frac{2\pi}{n}}\right\}\]

Hence the equation that solves the problem in this case,

\[1-\frac{l^2}{2r^2}=\cos\left(\frac{2\pi}{n}\right)\]

, whose positive result is

\begin{equation}
r=\frac{l}{\sqrt{2\left(1-\cos\left(\frac{2\pi}{n}\right)\right)}}=\frac{l}{2\sin\left(\frac{\pi}{n}\right)}
\label{eqn20.160}
\end{equation}

, that is known result that confirm again that method is correct.
\end{remark}

\section{Some ideas about the calculation of areas of convex polygons}

As it is known, the center of the circumscribed circle to a cyclic polygon can be found, depending on the lengths of the sides, within the area bounded by said sides (cyclic polygon centered on the inner (\emph{PCI})) or, on the contrary, outside of the area of the figure, or on its longest side (cyclic polygon not centered on the inner (\emph{PNCI}), was not found in the literature a classification for these cases). In the first case, for calculating the area , once calculated the radius of the circumscribed circle according to the method given here, the following sum must be applied,

\begin{equation}
A=\frac{1}{4} \sum_{k=1}^{n} l_k \sqrt{4r^2-l_k ^2 }
\label{eqn30,10}
\end{equation}

The terms of the sum are just the areas of triangles with sides $r$, $r$ and $l_k$, given in Figure 1. While in the second case you must subtract one of the terms, that is should be, 

\begin{equation}
A=\frac{1}{4} \sum_{k=1}^{n-1} l_k \sqrt{4r^2-l_k ^2 }-\frac{1}{4}l_n \sqrt{4r^2-l_n ^2 }
\label{eqn30,20}
\end{equation}

; where  $ l_n $ is the longer side. The latter is usually fairly obvious, however it should be noted that the longer side is, in the case corresponding to \emph{PNCI}, one that corresponds to the angle $\alpha_n=\sum^{n-1}_{1}\alpha_k $; it is not difficult to make sure that the triangle whose area is $\frac{l_{n}}{2}\sqrt{4r^{2}-l^{2}_{n}}$ is outside of the figure whose area is searched and yet his area is included in the first term of (\ref{eqn30,20}) and therefore must be subtracted from there.

\subsection{Some criteria for to determinate the position of the center of the circle with respect to convex polygon.}

\subsubsection{Necessary and sufficient criterion.} 
Considering expressions (\ref{eqn20,10}) and (\ref{eqn20.50}) we can easily infer the following necessary and sufficient condition so that the polygon is \emph{PNCI} (which derives from the necessity and sufficiency of (\ref{eqn20,10}) in this sense):

\begin{equation}
\arccos(1-\frac{l_{n}^2}{2r^{2}})=\sum^{n-1}_{k=1}\arccos(1-\frac{l_{k}^2}{2r^{2}})
\label{eqn30,30}
\end{equation}

; where $ r $ is the radius calculated in advance; and, considering angles exclusively in the range $\left[0,2\pi\right)$.

\begin{remark}
Note that, in the case of equilateral triangle,which it is obviously a \emph{PCI}, the above criteria leads to absurdity $1 = 2$.
\end{remark}

We can, to finalize, then define the following

\begin{criterion}
For a cyclic polygon it can be considered \emph{PNCI} is necessary and sufficient that may be made full (\ref{eqn30,30}) considering angles exclusively in the range $\left[0,2\pi\right)$.
\end{criterion}

\subsubsection{Necessary criterion.} 
From the analysis the following proposition is known:

\begin{lemma}
if the functions $f(r)$ and $g(r)$ have the same derivative they differ only in a constant.
\end{lemma}

Considering the relationship (\ref{eqn30,30}), and applying lemma $2$ we get the following:

\begin{equation}
\frac{d(\arccos(1-\frac{{l_{n}}^2}{2r^{2}}))}{dr}=\frac{d(\sum^{n-1}_{k=1}\arccos(1-\frac{{l_{k}}^2}{2r^{2}}))}{dr}
\end{equation}

And after simple transformations:

\[\frac{\frac{l_{n}^2}{r^{3}}}{\sqrt{1-(1-\frac{l_{n}^2}{2r^{2}})^{2}}}=\sum^{n-1}_{k=1}\frac{\frac{l_{k}^2}{r^{3}}}{\sqrt{1-(1-\frac{l_{k}^2}{2r^{2}})^{2}}}\]

\begin{equation}
\frac{1}{\sqrt{4(\frac{r}{l_{n}})^{2}-1}}=\sum^{n-1}_{k=1}\frac{1}{\sqrt{4(\frac{r}{l_{k}})^{2}-1}}
\label{eqn30,60}
\end{equation}

And (\ref{eqn30,60}) is expression that we was searching.

We can, then, define the following

\begin{criterion}
For a cyclic polygon it can be considered \emph{PNCI} is necessary that may be made full (\ref{eqn30,60}).
\end{criterion}

\subsubsection{Second necessary criterion.} 
From (\ref{eqn30,60}) it may be deducted other criteria necessary easier for classification as \emph{PNCI}:

To this end we transform (\ref{eqn30,60}) as follows,

\begin{equation}
1=\sum^{n-1}_{k=1}\frac{\sqrt{4(\frac{r}{l_{n}})^{2}-1}}{\sqrt{4(\frac{r}{l_{k}})^{2}-1}}
\label{eqn30,70}
\end{equation}

However

\begin{equation}
\sqrt{4\left(\frac{r}{l_{k}}\right)^{2}-1}=4\frac{A_{k}}{l_{k}^{2}}
\label{eqn30,80}
\end{equation}

, being $A_{k}$ the area of the triangle with sides $r$, $r$ and $l_{k}$ that can to see in the figure 1. Now, watching (\ref{eqn30,80}) we can see that all terms are positive in (\ref{eqn30,70}); so we can say the following criterion: 

\begin{criterion}
For a cyclic n-sided polygon is \emph{PNCI}, it is necessary that the following inequality is satisfied (Of this criterion it is deduced in an immediate way that any polygon with two equal sides, being them bigger or equal than the other ones, it will be \emph{PCI}):

\begin{equation}
\frac{\sqrt{4(\frac{r}{l_{n}})^{2}-1}}{\sqrt{4(\frac{r}{l_{k}})^{2}-1}}<1, \forall k<n
\end{equation}
\end{criterion}

\subsection{Integrated approach in the area calculation. Convex polygon}

The problem of calculating the area no difficulties once found the radius of the circle (We assume the radius $r$ calculated by the method described in the section 4 that can be seen below). This area can be calculated with relative ease using the tools of integral calculus. We believe that in this field the use of the comprehensive offers great advantages for the simplicity it provides in the formulas and the resulting possibilities of generalization to higher order polygons that generate the algebraic via virtually insurmountable difficulties.

Then we give some explanatory elements in relation to the above statements.

We consider the \emph{PCI} first. Figure 3 shows that, for each side of the polygon of n sides, can calculate the area of the shaded region. is not difficult to see that the same will,

\begin{equation}
A_{k}(r)=2\int^{l_{k}/2}_{0}\sqrt{r^{2}-x^{2}}dx-\frac{l_{k}}{2}\sqrt{4r^{2}-l_{k}^{2}}
\label{eqn30.100}
\end{equation}

And, after the integral calculation,

\begin{equation}
A_{k}(r)=r^{2}\left(1+\frac{l_{k}}{2}\sqrt{1-\frac{l_{k}^{2}}{4}}\right)\arcsin\left(\frac{l_{k}}{2}\right)-\frac{l_{k}}{2}\sqrt{4r^{2}-l_{k}^{2}}
\label{eqn30.110}
\end{equation}

Expression \ref{eqn30.100} can to be expressed by integrating by parts of most simple form as follow:

\begin{equation}
	A_{k}(r)=2\int^{l_{k}/2}_{0}\frac{x^{2}}{\sqrt{r^{2}-x^{2}}}dx
	\label{eqn30.120}
\end{equation}

The area of the polygon is then the complement of the sum of all $A_{i}(r)$:

\begin{equation}
	A(r)=\pi r^{2}-\sum^{n}_{k=1}A_{k}(r)
\label{eqn30.130}
\end{equation}

In the case of \emph{PNCI} happen to one side, the largest by the way (give it the designation $l_{n}$), requires no calculations of the shaded area of the figure, but its complement, ie,

\begin{equation}
A_{n}^{c}(r)=\pi r^{2}-A_{n}(r)
\end{equation}

Where $A_{n}$, is calculated according \ref{eqn30.110}. Hence the total area is definitely:

\begin{equation}
A(r)=\pi r^{2}-\sum^{n-1}_{k=1}A_{k}(r)-A_{n}^{c}(r)=2A_{n}(r)-\sum^{n}_{k=1}A_{k}(r)
\end{equation}

\begin{figure}[h]
\centering
\includegraphics[width=0.50\textwidth]{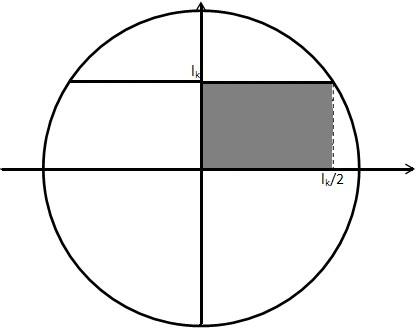}
\caption{Representing one side placed parallel to the $x$ axis.}
\label{fig:3}
\end{figure}

\section{Applications and examples}

\subsection{General formula for the triangle:}In this case the necessary expression for calculation is as fallows,

\[\left(1-\frac{a^2}{2r^2}+i\sqrt{1-\left(1-\frac{a^2}{2r^2}\right)^2 }\right)\left(1-\frac{b^2}{2r^2}+i\sqrt{1-\left(1-\frac{b^2}{2r^2}\right)^2 }\right)\times\]

\begin{center}
$\times$ $(1-\frac{c^{2}}{2 r^{2}}+i \sqrt{1-(1-\frac{c^2}{2r^{2}})^{2})}=1$
\end{center}

Taking the real part of the result we have to made the following equation,

\begin{center}
$-\sqrt{1 - (1 - \frac{a^2}{2 r^2})^2} \sqrt{1 - (1 - \frac{b^2}{2 r^2})^2} + \frac{a^2 b^2}{4 r^4} - \frac{a^2}{2 r^2} - \frac{b^2}{2 r^2} + \frac{c^2}{2 r^2}= 0$
\end{center}

And the solutions are as follows

\[r \rightarrow -\frac{a b c}{\sqrt{ 2 a^2 b^2+ 2 a^2 c^2 + 2 b^2 c^2 -a^4  - b^4 - c^4\}}}\]
\[r \rightarrow \frac{a b c}{\sqrt{2 a^2 c^2+ 2 b^2 c^2+ 2 a^2 b^2-a^4  - b^4   - c^4 }}\]

Taking the positive root, the result is consistent with existing theory. Factoring the denominator is the solution

\[r=\frac{abc}{\sqrt{(a+b+c)(b+c-a)(a+b-c)(a+c-b)}}\]

In the case of the equilateral triangle is the known result
\[r=\frac{a}{\sqrt{3}}\]

\subsubsection{Area of e triangle:}Having in main the last result we have the following

\[A=\frac{1}{4} a\sqrt{4 r^{2}-a^{2}}+\frac{1}{4} b\sqrt{4 r^{2}-b^{2}}+\frac{1}{4} c\sqrt{4 r^{2}-c^{2}}\]
\[=\frac{a^2 \sqrt{(-a^2 + b^2 + c^2)^2} + c^2 \sqrt{(a^2 + b^2 - c^2)^2} + b^2 \sqrt{a^2 - b^2 + c^2)^2}}{4\sqrt{(-a^4 + 2 a^2 b^2 - b^4 + 2 a^2 c^2 + 2 b^2 c^2 - c^4)}}=\]
\[=\frac{1}{4} \sqrt{(a + b - c) (a - b + c) (-a + b + c) (a + b + c)}\]

And if it do $s=\frac{a+b+c}{2}$ will be

\[A=\sqrt{(s - c) (s - b) ( s-a ) s}\]

And the both expression is coincident with the theory (Heron formula).

\subsection{General formula for the cyclic quadrilateral:}In this case the necessary expression for calculation is as fallows,

\[\left(1-\frac{a^2}{2r^2}+i\sqrt{1-\left(1-\frac{a^2}{2r^2}\right)^2 }\right)\times\left(1-\frac{b^2}{2r^2}+i\sqrt{1-\left(1-\frac{b^2}{2r^2}\right)^2 }\right)\times\]
\[\times\left(1-\frac{c^2}{2r^2}+i\sqrt{1-\left(1-\frac{c^2}{2r^2}\right)^2 }\right)\times\left(1-\frac{d^2}{2r^2}+i\sqrt{1-\left(1-\frac{d^2}{2r^2}\right)^2 }\right)=1\]

Taking the real part is the following equation:

\[-\sqrt{1 - (1 - \frac{a^2}{2 r^2})^2} \sqrt{1 - (1 - \frac{b^2}{2 r^2})^2} - \sqrt{1 - (1 - \frac{a^2}{2 r^2})^2} \sqrt{1 - (1 - \frac{c^2}{2 r^2})^2} -\]
\[-\sqrt{1 - (1 - \frac{b^2}{2 r^2})^2} \sqrt{1 - (1 - \frac{c^2}{2 r^2})^2} +\frac{c^2 \sqrt{1 - (1 - \frac{a^2}{2 r^2})^2} \sqrt{1 - (1 - \frac{b^2}{2 r^2})^2}}{2 r^2}+ \]
\[+ \frac{b^2 \sqrt{1 - (1 - \frac{a^2}{2 r^2})^2} \sqrt{1 - (1 - \frac{c^2}{2 r^2})^2}}{2 r^2}+ \frac{a^2 \sqrt{1 - (1 - \frac{b^2}{2 r^2})^2} \sqrt{1 - (1 - \frac{c^2}{2 r^2})^2}}{2 r^2}- \]
\[- \frac{a^2 b^2 c^2}{8 r^6} + \frac{a^2 b^2}{4 r^4} + \frac{a^2 c^2}{4 r^4} + \frac{b^2 c^2}{4 r^4} - \frac{a^2}{2 r^2} - \frac{b^2}{2 r^2} - \frac{c^2}{2 r^2} + \frac{d^2}{2 r^2} = 0\]

And the roots are as follows

\[r_{1}\rightarrow\frac{-\sqrt{-(bc+ad)(ac+bd)(ab+cd)}}{\sqrt{a-b-c-d)(a+b+c-d)(a+b-c+d)(a-b+c+d)}}\]
\[r_{2}\rightarrow\frac{\sqrt{-(bc+ad)(ac+bd)(ab+cd)}}{\sqrt{(a - b - c - d) (a + b + c - d) (a + b - c + d) (a - b + c + d)}}\]
\[r_{3} \rightarrow \frac{-\sqrt{(-b c + a d) (a c - b d) (a b - c d)}}{\sqrt{(a + b - c - d) (a - b + c - d) (a - b - c + d) (a + b + c + d)}}\]
\[r_{4} \rightarrow \frac{\sqrt{(-b c + a d) (a c - b d) (a b - c d)}}{\sqrt{(a + b - c - d) (a - b + c - d) (a - b - c + d) (a + b + c + d)}}\]

Performing a simple test in which $a=b=c=d=1$ we have

\[\left\{r_{1}\rightarrow -1.4142135623730950488\right\}\]
\[ \left\{r_{2} \rightarrow 0.7071067811865474617\right\}\]
\[\left\{r_{3} \rightarrow Indeterminate\right\}\]
\[\left\{r_{4} \rightarrow Indeterminate\right\}\]

From where $r_{2}$ is the solution, which it is consistent with the case of the square in that its sides it is all equal to 1; that is to say, the middle of diagonal line.
The quadrilateral and any regular polygon are \emph{PCI}, so its area is calculated by the expression (\ref{eqn30,10}). In this case the area is $1$, which can be determined easily.

\begin{remark}
From now on only specific cases will be presented, considering that the expression obtained for the pentagon $\left\{{a,a,a,b,c}\right\}$ it is already extremely complex and impossible to be analyzed explicitly with some practical purpose.
\end{remark}

\subsection{Cyclic Pentagon}

\subsubsection{Specific case for the cyclic pentagon of sides {1,2,4,5,5}:} In this case the necessary expression for calculation is as fallows,

\[\left(1-\frac{1^2}{2r^2}+i\sqrt{1-\left(1-\frac{1^2}{2r^2}\right)^2 }\right)\times\left(1-\frac{2^2}{2r^2}+i\sqrt{1-\left(1-\frac{2^2}{2r^2}\right)^2 }\right)\times\]
\[\times\left(1-\frac{4^2}{2r^2}+i\sqrt{1-\left(1-\frac{4^2}{2r^2}\right)^2 }\right)\times\left(1-\frac{5^2}{2r^2}+i\sqrt{1-\left(1-\frac{5^2}{2r^2}\right)^2 }\right)\times\]
\[\times\left(1-\frac{5^2}{2r^2}+i\sqrt{1-\left(1-\frac{5^2}{2r^2}\right)^2 }\right)=1\]

Taking the real part and similarly to the previous cases resolved is obtained the result

\[{r\rightarrow 3.04568}\]

\subsubsection{Calculation of area in square units:} Applying the criterion $3$ will be the following

\[\frac{\sqrt{4(\frac{r}{5})^{2}-1}}{\sqrt{4(\frac{r}{5})^{2}-1}}=1\]

But the criterion makes it clear that all these quotients, where the numerator contains the long side, should be less than unity; so this pentagon is \emph{PCI} (Anyone can to note from this moment that if a polygon has two equal sides and they are the longest then the polygon is \emph{PCI}. All regular polygon is \emph{PCI} of course). The expression to be applied, is then 
(\ref{eqn30,10}):

\[A=\frac{1}{4} \sum_{1}^{5} l_i \sqrt{4r^2-(l_i )^2 }=\]
\[=\frac{1}{4}u \sqrt{4 (3.0456755799776585)^2 - u^2} + \frac{1}{4} 2u\sqrt{4 (3.0456755799776585)^2 - \left(2u\right)^2}+\]
\[ + \frac{1}{4}4u \sqrt{4 (3.0456755799776585)^2 - \left(4u\right)^2}+ \frac{1}{4}5u \sqrt{4 (3.0456755799776585)^2 - \left(5u\right)^2}+\]
\[+\frac{1}{4}5u \sqrt{4 (3.0456755799776585)^2 - \left(5u\right)^2}\approx 17.6709 u^{2}\]

\subsubsection{Regular Pentagon} In this case the necessary expression for calculation is as fallows,

\[\left(1 - \frac{a^2}{2 r^2} + i \sqrt{1 - (1 - \frac{a^2}{2 r^2})^2}\right)^{5}=1\]

Taken the real part and solving 

\[\frac{a^8}{2 r^8} - \frac{4 a^6}{r^6} +\frac{10 a^4}{r^4} - \frac{15 a^2}{2 r^2} = 0\]

The above mentioned is equivalent to what continues,

\[a^2 (a^2 - 3 r^2) (a^4 - 5 a^2 r^2 + 5 r^4)= 0\]

And

\[r \rightarrow -a\frac{\sqrt{5  - \sqrt{5} }}{\sqrt{10}}, r \rightarrow a\frac{\sqrt{5  - \sqrt{5} }}{\sqrt{10}}\]
\[, r \rightarrow -a\frac{\sqrt{5  + \sqrt{5} }}{\sqrt{10}}, r \rightarrow a\frac{\sqrt{5  + \sqrt{5} }}{\sqrt{10}}\]

And the solution evidently is in the convex case,

\[r =\frac{\sqrt{5  + \sqrt{5} }}{\sqrt{10}}a\approx0.850651a\]

The solution would have been the same by applying the expression \ref{eqn20.160} and calculus although much simpler solution would have been lost for some non convex. But the other positive solution have sense to because are solution of non convex case, namely $r=\frac{\sqrt{5  - \sqrt{5} }}{\sqrt{10}}a\approx0.525731a$ is corresponding with fig.4.

\begin{figure}[h]
\centering
\includegraphics[width=0.30\textwidth]{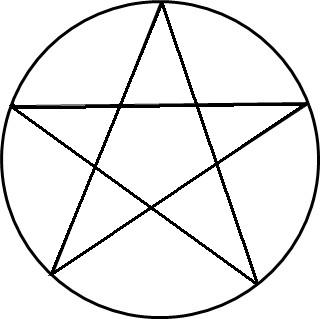}
\caption{Five pointed star.}
\label{fig:4}
\end{figure}

\subsection{Cyclic Heptagon:}

\subsubsection{Case to heptagon of sides $\left\{2,2,2,3,3,3,4\right\}$:} In this case the necessary expression for calculation is as fallows,

	\[\left(1-\frac{2^2}{2r^2}+i\sqrt{1-\left(1-\frac{2^2}{2r^2}\right)^2}\right)^{3}\left(1 - \frac{3^2}{2 r^2} + i \sqrt{1 - \left(1 - \frac{3^2}{2 r^2}\right)^2}\right)^{3}\times\]
	\[\times\left(1 - \frac{4^2}{2 r^2} + i \sqrt{1 - \left(1 - \frac{4^2}{2 r^2}\right)^2}\right)=1\]

Taking the real part and solving with respect to $r$ it shall be the roots

\[{\left\{r\rightarrow-3.15404\right\}}, {\left\{r \rightarrow 3.15404\right\}}\]

Obviously the root that is required is $r=3.15404$.

\subsubsection{Regular cyclic heptagon:} In this case the necessary expression for calculation is as fallows,

\[\left(1 - \frac{a^2}{2 r^2} + i \sqrt{1 - \left(1 - \frac{a^2}{2 r^2}\right)^2}\right)^{7} = 1\]

Taking the real part must indicate

\[\frac{a^{14}}{2 r^{14}} - \frac{7 a^{12}}{r^{12}} + \frac{77 a^{10}}{2 r^{10}} - \frac{105 a^{8}}{r^{8}} + \frac{147 a^6}{r^6} - \frac{98 a^4}{r^4} + \frac{49 a^2}{2 r^2}=0\]

Is obtained by factoring to following,

\[a^6 - 7 a^4 r^2 + 14 a^2 r^4 - 7 r^6= 0\]

This will shed extremely complicated solutions, however giving specific values to a, the result can be selected easily from the roots; as in the following case where $a=1$:

\[1 - 7 r^2 + 14 r^4 - 7 r^6 = 0\]

\[\left\{r \rightarrow -1.15238\right\}, \left\{r \rightarrow 1.15238 \right\}\]
\[ \left\{r \rightarrow -0.512858\right\},\left\{r \rightarrow 0.512858\right\}\]
\[\left\{r \rightarrow -0.639524\right\},\left\{r \rightarrow 0.639524\right\}\]

The only value that multiplied by $2\pi$ exceeds $7$, which is the perimeter of heptagon, is $1.15238$ giving concretely $1.1523824354812433*2\pi=7.240632386867575$, therefore this value is the desired result (of cols, in this case equation \ref{eqn20.160} can to be used for this purpose to).

Case $r \rightarrow 0.639524$ correspond to a non convex regular seven-pointed star that is represented in the figure 5,

\begin{figure}[hp]
\centering
\includegraphics[width=0.30\textwidth]{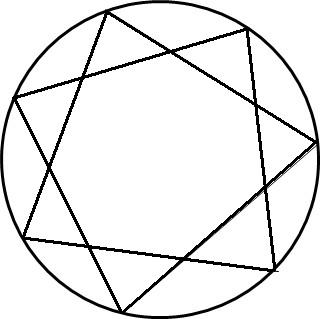}
\caption{Seven-pointed star}
\label{fig:5}
\end{figure}

To highlight the power of the method now consider some higher cases.

\subsection{Cyclic polygons of a superior quantity of sides:}

\subsubsection{Specific case for the polygon of 13 sides 1,1,1,1,3,3,3,3,3,3,3,4,6:} In this case the necessary expression for calculation is as fallows,

\[\left(1 - \frac{1^2}{2 r^2} + i \sqrt{1 - \left(1 - \frac{1^2}{2 r^2}\right)^2}\right)^{4} \times\left(1 - \frac{3^2}{2 r^2} + i \sqrt{1 - \left(1 - \frac{3^2}{2 r^2}\right)^2}\right)^{7}\times\]

\[\times\left(1 - \frac{4^2}{2 r^2} + i \sqrt{1 - \left(1 - \frac{4^2}{2 r^2}\right)^2}\right)\times\left(1 - \frac{6^2}{2 r^2} + i \sqrt{1 - \left(1 - \frac{6^2}{2 r^2}\right)^2}\right)=1\]

Taking the real part, the result is

\[\left\{r \rightarrow 5.67658\right\},\left\{ r \rightarrow -5.67658\right\},\]\[\left\{r \rightarrow 3.07337\right\},\left\{ r \rightarrow -3.07337\right\}\]

And the radio that we seek is clearly $r=5.676576550302839$ ($5.676576550302839*2\pi=35.667$ wile the perimeter indeed is 35 units). 

To calculate the area must determine if it is \emph{PCI} or \emph{PNCI}. We will use the criterion 1:

\[\arccos\left(1-\frac{l_{n}^2}{2r^{2}}\right)=\sum^{n-1}_{k=1}\arccos\left(1-\frac{l_{k}^2}{2r^{2}}\right)\]

The calculation gives the following results,

\[1.11364=5.16955\]

And the falsity of the previous equality analyzed indicates that the polygon is \emph{PCI}. Then we proceed using the expression (\ref{eqn30,10}):

\[A =\frac{1}{4}\sum^{13}_{k=1}l_{k}\sqrt{4r^{2}-l_{k}^{2}}\approx 93.8769 u^{2}\]

According to the integral expression (\ref{eqn30.130}), it is obtained, of course, the same result:

\[A =\pi r^{2}-\sum^{13}_{k=1}\left(2\int^{\frac{l_{k}}{2}}_{0}\sqrt{r^{2}-\xi^{2}}d\xi\right)+\sum^{13}_{k=1}\left(\frac{l_{k}}{2}\sqrt{4r^{2}-l_{k}^{2}}\right)\approx 93.8769 u^{2}\]

Both results have been placed here exactly like they were obtained in the package using, in their respective cases, the command 

\[\left\langle \left\langle ...\textnormal{ValuesData}\left[R(r,l_{1},l_{2},...l_{n}),l_{1}=v_{1},l_{2}=v_{2},...l_{n}=v_{n}\right]...\right\rangle\right\rangle\]

, where $R(r,l_{1},l_{2},...l_{n})$ represents, in this case, the area, dependent of the radius and the sides of the polygon.

\subsubsection{Regular polygon of 77 sides:} In this case the necessary expression for calculation is as fallows,

\[\left(1 - \frac{a^2}{2 r^2} + i\sqrt{ 1 - \left(1 - \frac{a^2}{2 r^2}\right)^2}\right)^{77}=1\]

The real part is, after factoring, the following:

\begin{center}
$a^2 (a^6 - 7 a^4 r^2 + 14 a^2 r^4 - 7 r^6)^2 (a^{10} - 11 a^8 r^2 + 44 a^6 r^4 - 77 a^4 r^6 + 55 a^2 r^8 - 11 r^{10})^2 (a^{60} - 59 a^{58} r^2 + 1652 a^{56} r^4 - 29205 a^{54} r^6 + 365859 a^{52} r^8 - 3455335 a^{50} r^{10} + 25556440 a^{48} r^{12} - 151794021 a^{46} r^{14} + 736647495 a^{44} r^{16} - 2956412711 a^{42} r^{18} + 9894941476 a^{40} r^{20} - 27773378270 a^{38} r^{22} + 65592924343 a^{36} r^{24} - 130530017729 a^{34} r^{26} + 218801812945 a^{32} r^{28} - 308337579027 a^{30} r^{30} + 363972489209 a^{28} r^{32} - 357976286928 a^{26} r^{34} + 291226209171 a^{24} r^{36} - 194130948828 a^{22} r^{38} + 104766369144 a^{20} r^{40} - 45083463663 a^{18} r^{42} + 15176439088 a^{16} r^{44} - 3900841911 a^{14} r^{46} + 742221909 a^{12} r^{48} - 100422335 a^{10} r^{50} + 9156888 a^8 r^{52} - 521938 a^6 r^{54} + 16580 a^4 r^{56} - 240 a^2 r^{58} + r^{60})^2$
\end{center}

Among all the roots are obtained, the only multiplied by $2\pi$ is higher than $77a$, is $12.2583 a$ ($r=77.02136692962065a$), so is the solution sought.

According to relation (\ref{eqn30.130}) the area is, in this case, $A\approx 471.551 a^2$ that which agrees almost perfectly with the area of the bounded circumference $\pi r^{2}\approx 472.074 a^2$.

\subsubsection{Regular polygon of 200 sides:}
This case that we won't reproduce here for space lack, shows as, for $200$ sides, the regular polygon can be considered, for the practical effects, similar to the circumscribed circle. In fact the obtained radius is $r=31.832297653000282 a$ being $a$ the longitude of the sides of the polygon, and the perimeter of the circle is $P=31.832297653000282 a*2 \pi=200.008 a$ differing, like it can be differentiated of the perimeter of the polygon in alone $0.008 a$.

\begin{remark}
In the two previous cases you can use the expression \ref{eqn20.160} but they would get lost the solutions for the non convex cases
\end{remark}

\subsection{Coinciding with example given in the literature:}
Ultimately, in $\left[1\right]$ pag.$228$ seven cyclic pentagons appear, one convex and six non-convex whose length of sides are 29, 30, 31, 32, y 33 respectively. Among the non-convex appear two irregular stars has five points, three with two crosses between sides, and one with a single crossing that resembles a fish. Of all radii of the respective circumscribed circles appear . The authors of this paper carried out the calculations of these radios by the method presented here resulting in a complete coincidence. 

The convex can be easily repeated case the exact steps performed in the examples above, however in the case of non-convex the expression must submit changes in given case. Here are all the cases occurring:

\subsubsection{Convex Case:} In this case the necessary expression for calculation is as fallows,

\[\left(1 - \frac{(29)^2}{2 r^2} + i \sqrt{1 - \left(1 - \frac{(29)^2}{2 r^2}\right)^2}\right) \left(1 - \frac{(30)^2}{2 r^2} + i \sqrt{1 - \left(1 - \frac{(30)^2}{2 r^2}\right)^{2}} \right) \times\]

\[\times\left(1 - \frac{(31)^2}{2 r^2} + i \sqrt{1 - \left(1 - \frac{(31)^2}{2 r^2}\right)^2}\right)\left(1 - \frac{(32)^2}{2 r^2} + i \sqrt{1 - \left(1 - \frac{(32)^2}{2 r^2}\right)^2}\right)\times\]

\[\times\left(1 - \frac{(33)^2}{2 r^2} + i \sqrt{1 - \left(1 - \frac{(33)^2}{2 r^2}\right)^2}\right)=1\]

Solutions: $r \rightarrow -26.38467157819376, r \rightarrow 26.38467157819376, r \rightarrow -16.512436619036198, r \rightarrow 16.512436619036198$.

Radius given by Robbins:

\begin{figure}[h]
\centering
\includegraphics[width=0.30\textwidth]{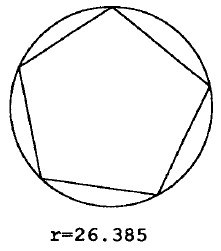}
\label{fig:6}
\end{figure}

\begin{remark}
As can be seen one solution is also coincident with the non-convex case of one of the stars. This is explained by the fact that, in the specific case of this star, the center of the circumcircle is in a position (with respect to polygon) such that angles in the exponent of (\ref{eqn20.130}) all 
have the same integer coefficient equal to $2$ (Each angle has the same number of sides passing through its own sector (See Lemma 1)), from here in (\ref{eqn20.150}) it is $E=2$ and the And the product for the calculation of r coincides. 
\end{remark}

\begin{figure}[h]
\centering
\includegraphics[width=0.20\textwidth]{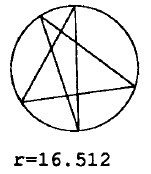}
\label{fig:7}
\end{figure}

\subsubsection{Non-Convex Cases:} Us consider the case in which only one of the angles has three sides in its sector of circumference. In that case it shall be, $3 \alpha_{1}+\alpha_{2}+\alpha_{3}+\alpha_{4}+\alpha_{5}=2\pi+2\alpha_{1}$; from here you can establish equality

\[e^{i\left(3 \alpha_{1}+\alpha_{2}+\alpha_{3}+\alpha_{4}+\alpha_{5}\right)}=e^{i2\alpha_{1}}\]

But as seen in (\ref{eqn20.130}) the above is equivalent to,

\[\prod^{n}_{p=1}N_{p}\left(r\right)=e^{i2\alpha_{1}}\]

Suppose further that $\alpha_{1}$ is such that $e^{i\alpha_{1}}=N_{1}\left(r\right)$, you are then obtained the formula,

\[\prod^{n}_{p=1}N_{p}\left(r\right)=N_{1}^{2}\left(r\right)\]

Now multiplied by the complex conjugate $N_{1}^{*}\left(r\right)$ is that,

\begin{equation}
N_{1}^{*}\left(r\right)\prod^{n}_{p=2}N_{p}\left(r\right)=1
\label{eqn40.10}
\end{equation}

Implemented with (\ref{eqn40.10}) a process completely analogous to (\ref{eqn20.30}) and changing $N_{1}^{*}\left(r\right)$ as different sides have corroborated all results Robbins as follows:

\begin{itemize}

\item ------------------------------------------------------------
	\[N_{1}^{*}\left(r\right)=1-\frac{(29)^{2}}{2 r^2}- i \sqrt{1 - (1-\frac{(29)^{2}}{2 r^2})^2}\]

Our results: $r=18.651420360146222$.

Robbins result:

\begin{figure}[h]
\centering
\includegraphics[width=0.20\textwidth]{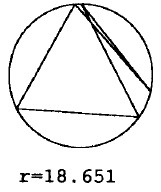}
\label{fig:8}
\end{figure}

\item --------------------------------------------------------------
	\[N_{1}^{*}\left(r\right)=1-\frac{(30)^{2}}{2 r^2}- i \sqrt{1 - (1-\frac{(30)^{2}}{2 r^2})^2}\]

Our results: $r=18.33473396357385$.

Robbins result:

\begin{figure}[h]
\centering
\includegraphics[width=0.20\textwidth]{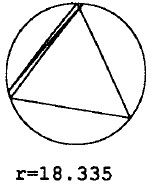}
\label{fig:8}
\end{figure}

\item --------------------------------------------------------------
	\[N_{1}^{*}\left(r\right)=1-\frac{(31)^{2}}{2 r^2}- i \sqrt{1 - (1-\frac{(31)^{2}}{2 r^2})^2}\]

Our results: $r=17.99086611423847$.

Robbins result:

\begin{figure}[h]
\centering
\includegraphics[width=0.20\textwidth]{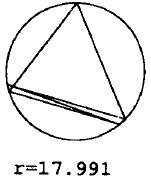}
\label{fig:9}
\end{figure}

\item --------------------------------------------------------------
	\[N_{1}^{*}\left(r\right)=1-\frac{(32)^{2}}{2 r^2}- i \sqrt{1 - (1-\frac{(32)^{2}}{2 r^2})^2}\]

Our results: $r=17.595145969431748$.

Robbins result:

\begin{figure}[h]
\centering
\includegraphics[width=0.20\textwidth]{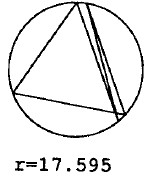}
\label{fig:10}
\end{figure}

\item --------------------------------------------------------------
	\[N_{1}^{*}\left(r\right)=1-\frac{(33)^{2}}{2 r^2}- i \sqrt{1 - (1-\frac{(33)^{2}}{2 r^2})^2}\]

Our results: $r=17.02586455404377$.

Robbins result:

\begin{figure}[h]
\centering
\includegraphics[width=0.20\textwidth]{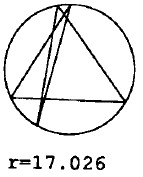}
\label{fig:11}
\end{figure}

\end{itemize}

\section{Conclusions}

\begin{conclusion}
It has presented a method that allows, in principle, find the radius of the circle corresponding to any convex cyclical polygon, as well as consequently the area of said figure.
\end{conclusion}

\begin{conclusion}
The presented method, in its most general expression allows also find the radius of the circumscribed circumferences to all non convex polygon cyclical.
\end{conclusion}

\begin{conclusion}
Some necessary and sufficient criterions have been determined for the determination of the position of the center of the circumference bounded to the cyclic polygon, that which is important for the calculation of the areas and in next works they will be key to develop a more efficient similar method.
\end{conclusion}

\begin{conclusion}
The general expressions for radius on $n>4$ tend to be of a  degree of complexity prohibitive of their use  in explicit form, but the calculation of these values for specific polygons it not represents no difficulty with the use of this method, unless the hardware processing capacity is insufficient.
\end{conclusion}

\end{document}